\documentclass[]{article}
\usepackage[english]{babel}  
\usepackage{lineno} 
\usepackage{float}
\usepackage{graphicx,multicol} 
\usepackage{epic,eepic,epsfig} 
\usepackage{amssymb} 
\usepackage{tikz}
\usepackage{lineno}
\newcommand{\qed}{\hfill$\diamond$\\\vspace{2mm}}
\newcommand{\pf}{{\bf Proof: }} 
\newcommand{\be}{\begin{enumerate}} 
\newcommand{\ee}{\end{enumerate}} 
\newcommand{\bd}{\begin{description}} 
\newcommand{\ed}{\end{description}}

\newcommand{\beq}{\begin{equation}} 
\newcommand{\eeq}{\end{equation}}

\newtheorem{theorem}{Theorem}[section] 
\newtheorem{lemma}[theorem]{Lemma}

\newtheorem{corollary}[theorem]{Corollary}

\newtheorem{problem}[theorem]{Problem}

\newcommand{\AY}[1]{{#1}}

\begin{document} 
 \bibliographystyle{plain}
%\setpagewiselinenumbers 
%\modulolinenumbers[1] 
%\linenumbers 

 \title{Acyclic orientations of mixed graphs}
 \author{J\o{}rgen Bang-Jensen\thanks{Department of Mathematics and Computer Science, University of Southern Denmark, Odense Denmark and School of Mathematics, Shandong University, Jinan, China (email:jbj@imada.sdu.dk)} and Anders Yeo\thanks{Department of Mathematics and Computer Science, University of Southern Denmark, Odense Denmark and Department of Mathematics, University of Johannesburg, South Africa (email:yeo@imada.sdu.dk)}}
 \maketitle

\begin{abstract}
    A mixed graph $M=(V,E\cup A)$ is acyclic if its directed part $(V,A)$ is an acyclic digraph. In this note we study the so-called orientation completion problem for the class of acyclic mixed graphs. That is, given an acyclic mixed graph $M$ and a property ${\cal P}$; can we orient the edges of $M$ so that the resulting digraph is acyclic and has property ${\cal P}$. We prove that one can decide in polynomial time whether $M$ can be completed to an acyclic digraph with an out-branching from a prescibed vertex $s$, while it is NP-complete to decide whether $M$ has an acyclic orientation which has both an out-branching and an in-branching (a bipolar orientation). We  show that it is NP-complete to decide whether $M$ can be oriented so that it contains a directed path between two prescribed vertices. Finally we describe  a polynomial algorithm for deciding whether an acyclic digraph $D$ has an out-branching $B^+_s$ such that the digraph $D-A(B^+_s)$ is connected (in the underlying sense). Based on this we pose as an open  problem the complexity of deciding whether the edges of an acyclic mixed graph can be oriented so that the result is an  acyclic digraph with a non-separating out-branching.
    
\end{abstract}

\section{Introduction}

All notation and terminology used is consistent with \cite{bang2009}. Orientations of mixed graphs so as to satisfy prescribed reachability conditions or connectivity properties has been studied many times in the literature, see e.g. \cite{arkinDAM116,boeschAMM87,elberfeldTCS483, gamzuA74,hoerscharxiv2025,hoerschDO48}.
 Several of those problems fit the framework which is called the  {\bf orientation completion problem} introduced by 
 Bang-Jensen, Huang and Zhu in \cite{bangJGT87}. This is defined as follows: given 
 a mixed graph $M=(V,E\cup A)$ (which can be seen as a partially oriented graph) and a graph property ${\mathcal P}$; can we orient the edges of $E$ as an arc set $A'$ so that the resulting directed graph $D=(V,A\cup A')$ has property ${\mathcal P}$? The orientation completion problem framework generalizes many problems dealing with either orientations of undirected graphs or finding prescribed structures in digraphs. Thus the orientation completion problem for a given property ${\mathcal P}$ is at least as hard as both the orientation problem where one seeks an orientation of an undirected graph so that the \AY{resulting digraph} has property ${\mathcal P}$ and the already oriented version where we wish to know whether a given digraph has property ${\mathcal P}$. 
 \\
 
 The purpose of this note is to study some natural orientation completion problems when the property ${\mathcal P}$ includes that the resulting digraph $D$ must be acyclic.  In Section \ref{sec:outbr} \AY{we} first study whether it is possible to obtain a completion to an  acyclic digraph so that  a prescribed vertex $s$ has a directed path to all other vertices in $D$. Then we show in Section \ref{sec:bipolar} that it is NP-complete to decide whether there is a \AY{completion} to an acyclic digraph which has both a vertex that can reach other vertices by directed paths and another vertex that can be reached from all other vertices. In Section \ref{sec:stpath} we show that when we just want to ensure reachability between a specific pair of vertices (rather than having an out-branching), then the orientation completion problem to an acyclic digraph becomes NP-complete. Finally, we show in Section \ref{sec:nonsepBr} that one can decide in polynomial time whether an acyclic digraph has a non-separating out-branching and \AY{we, as} an interesting open problem, \AY{suggest} to study the complexity of \AY{the} orientation completion version of the problem.

\section{Completing an acyclic mixed graph to an acyclic digraph with an out-branching}\label{sec:outbr}

%The results in this section arose out of discussions with Jing Huang.\\

Let $M=(V,A,E)$ be an acyclic mixed graph and let $s\in V$ be given. A 3-partition $(Y,B,X)$ of $V$ is an {\bf $s$-obstruction} if 
\begin{enumerate}
  \item[(i)] $s\in Y\cup X$.
\item[(ii)] There are no edges or arcs with one end in $Y$ and the other in $B$.
  \item[(iii)] Every vertex in $X$ is reachable by a directed path from some vertex in $B$.
\end{enumerate}

\begin{lemma}
  \label{obstruction}
  \AY{Let} $M=(V,A,E)$ be an acyclic mixed graph and let $s\in V$ be given. If $M$ has an $s$-obstruction, then $M$ has no completion to an acyclic digraph with an out-branching rooted at $s$.
\end{lemma}
\pf Let $(Y,B,X)$ be an $s$-obstruction and let $B'\subseteq B$ denote the set of vertices from $B$ that can reach a vertex in $X$. Recall that, by (iii), every vertex in $X$ is reachable from $B'$. Let  $D$ be an  arbitrary acyclic completion of $M$ and let ${\cal O}=v_1,v_2,\ldots{},v_n$ be any acyclic ordering of $D$. Then ${\cal O}$ induces an ordering $b_1,b_2,\ldots{},b_k$, $k=|B'|$ of the vertices in $B'$. Since every vertex of $X$ is reachable from $B'$ in $M$, all vertices of $X$ must occur after $v_i=b_1$ in $\cal O$. This implies that there is no path from $s$ to $b_1$ in $D$ so $D$ has no out-branching from $s$. \qed

\begin{lemma}
  \label{extendor}
  If the directed part of an acyclic mixed graph $M=(V,A,E)$ has an out-branching $B^+_s$ rooted at $s$, then $M$ has a completion to an acyclic digraph containing $B^+_s$.
\end{lemma}
\pf
Let ${\cal O}=v_1,v_2,\ldots{},v_n$ be an acyclic ordering of $B^+_s$. For every edge $v_iv_j$ of $E$ we orient the edge from $v_i$ to $v_j$ if $i<j$ and from $v_j$ to $v_i$ if $i>j$. Clearly the resulting digraph is acyclic and it contains $B^+_s$. \qed

\begin{theorem} \label{thmObs}
  Let $M=(V,A,E)$ be an acyclic mixed graph and let $s\in V$ be given. $M$ can be oriented as an acyclic digraph with an out-branching rooted at $s$ if and only if $M$ has no $s$-obstruction.
\end{theorem}

\pf By Lemma \ref{obstruction} we may assume that $M$ has no $s$-obstruction.
This implies that no arc of $A$ enters $s$, as otherwise $(\emptyset{},V-s,\{s\})$ would be an $s$-obstruction. The following algorithm constructs a completion of $M$ to an acyclic digraph with an out-branching from $s$.

Let $Z$ denote the set of vertices that are reachable from $s$ via directed paths in $M$. If $Z=V$ we are done by Lemma \ref{extendor}. So assume $Z\neq V$. Let $Z'\subseteq Z$ denote the set of vertices in $Z$ \AY{which are adjacent to a vertex in $V-Z$ (by either an arc or an edge). Note that by the definition of $Z$ there are no arcs from $Z'$ to $V-Z$.}
If every vertex of $Z'$ is reachable from $V-Z$ then $(Z-Z',V-Z,Z')$ is an $s$-obstacle, contradicting the assumption. Hence some vertex $z\in Z'$ cannot be reached from $V-Z$. Let $zv$ be an arbitrary edge from $z$ to $V-Z$ and orient this edge from $z$ to $v$. Since $z$ was not reachable from $V-Z$, the directed part of the current mixed graph is acyclic. Now update $Z$ to include $v$ and all other vertices of $V-Z$ that became reachable when we oriented $zv$ from $z$ to $v$ and repeat the argument above.

\AY{Let $Z_1$ denote the updated $Z$ and $Z_1'$ be the updated $Z'$ and let $M_1$ be the updated $M$ (where the edge $zv$ has been oriented).  Note that if there is a directed path, in $M_1$, from a vertex in $V-Z_1$ to a vertex in $Z_1'$, then this is also a directed path in $M$ as if the arc $zv$ is on the path in $M_1$ then there is a path from $V-Z_1$ to $z$ implying that $z \not\in Z'$, a contradiction. Therefore, if every vertex of $Z'$ is reachable from $V-Z$ then $(Z-Z',V-Z,Z')$ is an $s$-obstacle also in later iterations of the algorithm.}

As at least one new vertex is added in each round, after at most $n-1$ steps, this results in an acyclic  mixed graph $M'$ with an out-branching rooted at $s$ and then we can complete the orientation by Lemma \ref{extendor}.\qed

\AY{We note that the proof of Theorem~\ref{thmObs} implies a polynomial algorithm for deciding if a mixed graph has an acyclic completion containing an out-branching rooted at $s$.}

\section{Completing to a bipolar orientation}\label{sec:bipolar}

A vertex ordering $v_1,v_2,\ldots{},v_n$ of a graph $G=(V,E)$ with two distinct vertices $s,t$ is called an 
{\bf $st$-ordering} of $V$ if it satisfies that $v_1=s,v_n=t$ and for every $1<j<n$ there are indices $s\leq i<j<k\leq n$ such that $v_iv_j\in E$ and $v_jv_k\in E$.

\begin{theorem}\label{stexists}\cite{lempel1967}
  \begin{itemize}
  \item Every 2-connected graph $G$ has an $st$-ordering for every choice of distinct vertices $s,t$ of $G$.
    \item A graph $G=(V,E)$ has an $st$-ordering if and only if the graph $G+st$ is 2-edge-connected.
    \end{itemize}
\end{theorem}

A {\bf bipolar orientation} of a graph $G=(V,E)$ is an acyclic orientation $D=(V,A)$ of $G$ which has exactly one source $s$ and one sink $t$. Thus, for every acyclic ordering $s=v_1,v_2,\ldots{},v_n=t$ of $D$ each $v_i$ with $1<i<n$ has an in-neighbour $v_j$ (which will have $j<i$)  and an out-neighbour $v_k$ (which will have $k>i$). From this it is easy to see that $G=(V,E)$ has a bipolar orientation with source $s$ and sink $t$ if and only if $G$ has an $st$-ordering.

\begin{theorem}
  It is NP-complete to decide whether an acyclic mixed graph $M=(V,E\cup A)$ can be completed to a bipolar orientation with source $s$ and sink $t$ for prescribed vertices $s,t$.
\end{theorem}

\pf We will reduce  the so-called {\bf betweenness} problem to our problem. The betweenness problem is as follows: given a collection ${\cal C}$ of tripels $(v_{i_1},v_{j_1},v_{k_1}),\ldots{},(v_{i_p},v_{j_p},v_{k_p})$ so that all elements of the tripels are from the set $W$; is it possible to order the elements of $W$ as ${\cal O}= v_1,v_2,\ldots{},v_n$ so that for every tripel $(v_{i_q},v_{j_q},v_{k_q})$, $q\in [p]$, the elements occur in ${\cal O}$ either in the order $v_{i_q},v_{j_q},v_{k_q}$ or in the order $v_{k_q},v_{j_q},v_{i_q}$?
This problem, listed as problem  [MS1] in \cite{garey1979}, is NP-complete.

Given an instance \AY{$(W,{\cal C})$} of betweenness where $W=\{v_1,v_2,\ldots{},v_n\}$ we construct a mixed acyclic graph $M=(V^*,E\cup A)$ where $E$ and $A$ where  we construct $V^*,E,A$ as follows starting from  $V^*=W\cup \{s,t\}$, $E=\emptyset{}$ and $A=\emptyset$.
\begin{itemize}
\item for every $v_i\in W$ add the arcs $sv_i$ and $v_it$ to $A$
  \item for every triple $(v_{i_q},v_{j_q},v_{k_q})$ of ${\cal C}$ we add two new vertices $x_q,y_q$ to $V^*$, add the four edges $v_{i_q}x_q,v_{i_q}y_q,x_qv_{k_q},y_qv_{k_q}$ to $E$ and the arcs $x_qv_{j_q},v_{j_q}y_q$ to $A$. Let $M_q$ denote the  mixed subgraph of $M$ induced by the vertices $x_q,y_q,v_{i_q},v_{j_q},v_{k_q}$.
\end{itemize}

It is easy to check that the resulting mixed graph $M$ has no directed cycle. Clearly $s$ will be a source and $t$ a sink in every orientation of $M$ (as an acyclic digraph). It remains to show that $M$ can be completed to an acyclic digraph with exactly one source and one sink if and only if \AY{$(W,{\cal C})$} is a 'yes'-instance of betweenness. Suppose first that $M$ has a bipolar orientation $D$. Since $M$ contains the directed path $x_qv_{j_q}y_q$ for every triple $(v_{i_q},v_{j_q},v_{k_q})$ of ${\cal C}$ and $x_q,y_q$ are private to $M_q$, the only possible acyclic orientations of $M_q$ are obtained by orienting the 4 edges as $v_{i_q}x_q, v_{i_q}y_q,x_qv_{k_q},y_qv_{k_q}$ or as
$v_{k_q}x_q, v_{k_q}y_q,x_qv_{i_q},y_qv_{i_q}$. In the first case the vertices will occur in the order $v_{i_q},x_q,v_{j_q},y_q,v_{k_q}$ in every acyclic ordering of $V^*$ and in the second case they will appear in the order $v_{k_q},x_q,v_{j_q},y_q,v_{i_q}$ in every acyclic ordering of $V^*$. Hence if $D$ has a bipolar orientation, then we obtain
the desired ordering ${\cal O}$ of $W$ by taking any acyclic order of $V^*$ and deleting all the vertices of $\{x_q|q\in [p]\}\cup\{y_q|q\in [p]\}$.

Suppose now that ${\cal O}$ is an ordering of $W$ that certifies that \AY{$(W,{\cal C})$} is a 'yes'-instance of betweenness, then we obtain a bipolar orientation of $M$ by orienting the four edges of $M_q$ as $v_{i_q}x_q, v_{i_q}y_q,x_qv_{k_q},y_qv_{k_q}$ if the vertices of $(v_{i_q},v_{j_q},v_{k_q})$ occur in the order $v_{i_q}<v_{j_q}<v_{k_q}$ in ${\cal O}$  and as
$v_{k_q}x_q, v_{k_q}y_q,x_qv_{i_q},y_qv_{i_q}$ if they occur in the order $v_{k_q}<v_{j_q}<v_{i_q}$ in $\cal O$. \qed

\section{Orienting a acyclic mixed graph as an acyclic digraph with an $(s,t)$-path}\label{sec:stpath}

We consider the following problem.

\begin{problem}[{\sc Completion to  an acyclic digraph with an $(s,t)$-path}]
  Given a mixed graph $M$ and two vertices $s,t$. Can $M$ be completed to an acyclic digraph $D$ with an $(s,t)$-path?
  \end{problem}

\medskip

%The following result arose out of a communication with Mathias Kriesell.

\begin{theorem}\label{thm:mixedpath}
{\sc Completion to an acyclic digraph with an $(s,t)$-path.} is NP-complete.
\end{theorem}

\pf Let $L$ be the mixed graph on 4 vertices $u,v,w,z$, two edges $uv,wz$ and two arcs
$(v,w)$, $(z,u)$. There are 4 possible completions of $L$ as an oriented graph and all except the one where we orient as $(u,v),(w,z)$ is acyclic. Let $W$ denote the mixed graph on 8 vertices $p,q,a_1,b_1,a_2,b_2,a_3,b_3$, three edges $a_1b_1,a_2b_2,a_3b_3$ and six arcs $(p,a_1),(p,a_2),(p,a_3),(b_1,q),(b_2,q),(b_3,q)$.

Let $\cal F$ be an instance of 3-SAT with variables $x_1,x_2,\ldots{},x_n$ and clauses $C_1,C_2,\ldots{},C_m$, where each clause is of the form $(\ell_1\vee\ell_2\vee\ell_3)$ and each $\ell_i$ is either one of the variables $x_j$ or the negation $\bar{x}_j$ of such a variable.

Let $p_i$ ($q_i$) be the number of times variable $x_i$ ($\bar{x}_i$) 
occurs as a literal in $\cal F$. 
The enumeration of the clauses $C_1,\ldots{},C_m$ induces an ordering on 
the occurrences of the same literal in the formula.
Guided by this ordering we  now construct a mixed graph $H'=H'({\cal F})$ as follows: For each variable $x_i$, $i\in [n]$ let $X_i$ denote the mixed graph that we obtain from taking the union of the two internally disjoint $(\alpha_i,\beta_i)$-paths $\alpha_iu_{i,1}v_{i,1}u_{i,2}v_{i,2}\ldots{}u_{i,p_i}v_{i,p_i}\beta_i$ and $\alpha_i\bar{u}_{i,1}\bar{v}_{i,1}\bar{u}_{i,2}\bar{v}_{i,2}\ldots{}\bar{u}_{i,q_i}\bar{v}_{i,q_i}\beta_i$ and then replacing each of the arcs of the kind $u_{i,j}v_{i,j}$ or $\bar{u}_{i,j}\bar{v}_{i,j}$ by an edge. For each clause $C_j$, $j\in [m]$ let $W_j$ denote a copy of $W$ where we denote  the vertices as follows 
$p_j,q_j,a_{j,1},b_{j,1},a_{j,2},b_{j,2},a_{j,3},b_{j,3}$. Now we form $H'$ from $X_1,\ldots{},X_n, W_1,\ldots{},W_m$ as follows
\begin{itemize}
\item add the arcs $\beta_i\beta_{i+1}$ for $i\in [n-1]$
\item add the arcs $q_jp_{j+1}$ for $j\in [m-1]$.
\item For each clause $C_j=(\ell_{j,1}\vee\ell_{j,2}\vee\ell_{j,3})$, $j\in [m]$ we add the 2 directed arcs of three copies of $L$ as follows (for $k\in [3]$): If $\ell_{j,k}=x_i$ and this is the $g$'th occurence of $x_i$ according to the ordering above then we add the arcs $(v_{i,g},a_{j,k}),(b_{j,k},u_{i,g})$. If $\ell_{j,k}=\bar{x}_f$ this is the $h$'th occurence of $x_f$ according to the ordering above then we add the arcs $(\bar{v}_{f,h},a_{j,k}),(b_{j,k},\bar{u}_{f,h})$.
\item Finally we add the arc $(\beta_n,p_1)$ and let $s=\alpha_1,t=q_m$.
\end{itemize}

\begin{figure}[H]
\begin{center}

\tikzstyle{vertexB}=[circle,draw, minimum size=15pt, scale=0.6, inner sep=0.5pt]
\begin{tikzpicture}[scale=0.55,>=stealth]
\node (s) at (-1,10) [vertexB] {$s$};
\node (a1) at (0.5,12) [vertexB] {};
\node (b1) at (3,12) [vertexB] {};
\node (a'1) at (0.5,8) [vertexB] {};
\node (b'1) at (3,8) [vertexB] {};
\node (u1) at (4.5,10)[vertexB] {};
\node (v1) at (7,10) [vertexB] {};
\node (a2) at (8.5,12) [vertexB] {};
\node (b2) at (11,12) [vertexB] {};
\node (a'2) at (8.5,8) [vertexB] {};
\node (b'2) at (11,8) [vertexB] {};
\node (u2) at (12.5,10)[vertexB] {};
\node (v2) at (15,10) [vertexB] {};
\node (a3) at (16.5,12) [vertexB] {};
\node (b3) at (19,12) [vertexB] {};
\node (a'3) at (16.5,8) [vertexB] {};
\node (b'3) at (19,8) [vertexB] {};
\node (u3) at (20.5,10)[vertexB] {};

\draw [->] (s) to (a1);
\draw [->] (s) to (a'1);
\draw (a1) to (b1);
\draw [->] (b1) to (u1);
\draw (a'1) to (b'1);
\draw [->] (b'1) to (u1);
\draw [->] (u1) to (v1);
\draw [->] (v1) to (a2);
\draw [->] (v1) to (a'2);
\draw (a2) to (b2);
\draw [->] (b2) to (u2);
\draw (a'2) to (b'2);
\draw [->] (b'2) to (u2);
\draw [->] (u2) to (v2);
\draw [->] (v2) to (a3);
\draw [->] (v2) to (a'3);
\draw (a3) to (b3);
\draw [->] (b3) to (u3);
\draw (a'3) to (b'3);
\draw [->] (b'3) to (u3);

\node (t) at (0,0) [vertexB] {$t$};
\node (d3) at (1.5,2) [vertexB] {};
\node (d'3) at (1.5,0) [vertexB] {};
\node (d''3) at (1.5,-2) [vertexB] {};
\node (c3) at (4,2) [vertexB] {};
\node (c'3) at (4,0) [vertexB] {};
\node (c''3) at (4,-2) [vertexB] {};
\node (w3) at (5.5,0) [vertexB] {};
\node (z2) at (7,0) [vertexB] {};
\node (d2) at (8.5,2) [vertexB] {};
\node (d'2) at (8.5,0) [vertexB] {};
\node (d''2) at (8.5,-2) [vertexB] {};
\node (c2) at (11,2) [vertexB] {};
\node (c'2) at (11,0) [vertexB] {};
\node (c''2) at (11,-2) [vertexB] {};
\node (w2) at (12.5,0) [vertexB] {};
\node (z1) at (14,0) [vertexB] {};
\node (d1) at (15.5,2) [vertexB] {};
\node (d'1) at (15.5,0) [vertexB] {};
\node (d''1) at (15.5,-2) [vertexB] {};
\node (c1) at (18,2) [vertexB] {};
\node (c'1) at (18,0) [vertexB] {};
\node (c''1) at (18,-2) [vertexB] {};
\node (w1) at (19.5,0) [vertexB] {};

\draw [->] (w1) to (c1);
\draw [->] (w1) to (c'1);
\draw [->] (w1) to (c''1);
\draw (c1) to (d1);
\draw (c'1) to (d'1);
\draw (c''1) to (d''1);
\draw [->] (d1) to (z1);
\draw [->] (d'1) to (z1);
\draw [->] (d''1) to (z1);
\draw [->] (z1) to (w2);
\draw [->] (w2) to (c2);
\draw [->] (w2) to (c'2);
\draw [->] (w2) to (c''2);
\draw (c2) to (d2);
\draw (c'2) to (d'2);
\draw (c''2) to (d''2);
\draw [->] (d2) to (z2);
\draw [->] (d'2) to (z2);
\draw [->] (d''2) to (z2);
\draw [->] (z2) to (w3);
\draw [->] (w3) to (c3);
\draw [->] (w3) to (c'3);
\draw [->] (w3) to (c''3);
\draw (c3) to (d3);
\draw (c'3) to (d'3);
\draw (c''3) to (d''3);
\draw [->] (d3) to (t);
\draw [->] (d'3) to (t);
\draw [->] (d''3) to (t);
\draw [->] (u3) to (w1);

% The coloured arcs are straight unless a small displacement is needed
% to clear an intermediate vertex.
\draw [red,->] (b1) .. controls (8.04,8.73) and (13.04,5.40) .. (c1);
\draw [red,->] (d1) to (a1);
\draw [red,->] (b2) .. controls (13.26,7.96) and (15.60,3.96) .. (c'1);
\draw [red,->] (d'1) .. controls (12.92,3.86) and (10.59,7.86) .. (a2);
\draw [red,->] (b3) .. controls (19.66,7.26) and (19.33,2.60) .. (c''1);
\draw [red,->] (d''1) .. controls (16.83,2.60) and (17.16,7.26) .. (a3);

\draw [blue,->] (b1) .. controls (6.07,8.99) and (8.74,5.66) .. (c2);
\draw [blue,->] (d2) to (a1);
\draw [blue,->] (b'2) .. controls (11.58,5.33) and (11.58,2.67) .. (c'2);
\draw [blue,->] (d'2) .. controls (9.08,2.67) and (9.08,5.33) .. (a'2);
\draw [blue,->] (b'3) .. controls (16.74,4.34) and (14.07,1.01) .. (c''2);
\draw [blue,->] (d''2) to (a'3);

\draw [green,->] (b'1) to (c3);
\draw [green,->] (d3) to (a'1);
\draw [green,->] (b'2) to (c'3);
\draw [green,->] (d'3) to (a'2);
\draw [green,->] (b3) .. controls (13.81,7.54) and (8.81,2.87) .. (c''3);
\draw [green,->] (d''3) .. controls (6.36,2.81) and (11.36,7.48) .. (a3);
\end{tikzpicture}

\end{center}
\caption{The digraph $H'=H'({\cal F})$ when $\cal F$ consists of the clauses
$(x_1\vee{}x_2\vee{}x_3), (x_1\vee{}\bar{x}_2\vee{}\bar{x}_3), (\bar{x}_1\vee{}\bar{x}_2\vee{}x_3)$.}\label{stmixed graph}
\end{figure}
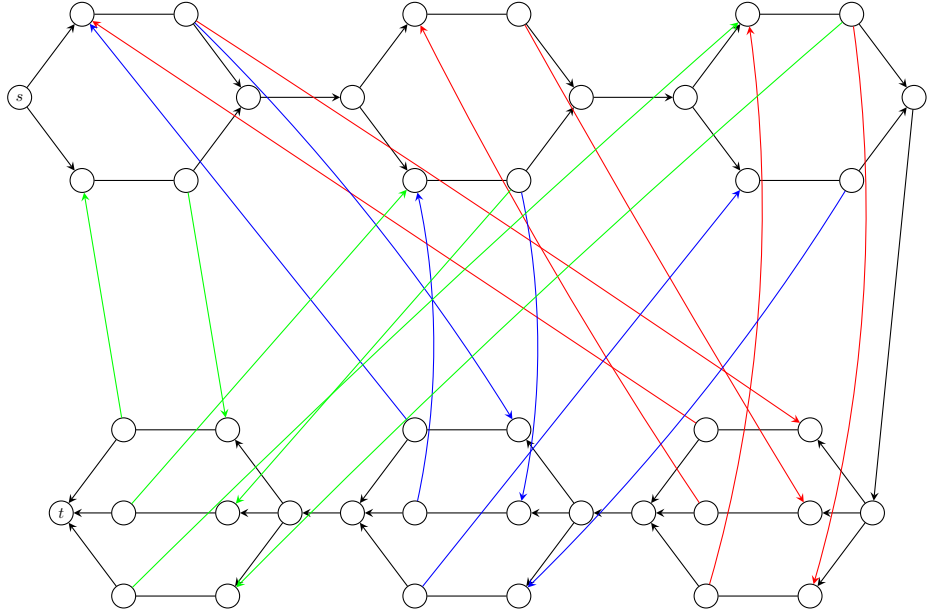

See Figure \ref{stmixed graph}.

%\end{figure}

We claim that $\cal F$ is satisfiable if and only if $H'$ has a completion to an acyclic digraph with an $(s,t)$-path. Suppose first that $H'$ has such a completion $D$ and let $P$ be an $(s,t)$-path in $D$. Then, for each $j\in [m]$ $D$ contains at least one of the arcs $(a_{j,h},b_{j,h})$. Fix such an arc $(a_{j,h(j)},b_{j,h(j)})$ for each $j\in [m]$ By the property of $L$ the arc corresponding to the litteral $\ell_{j,h(j)}$ is oriented backwards (towards the $u$ or $\bar{u}$ vertex), implying that it cannot be part of any $(s,t)$ path. This means that $P$ avoids at least one of the three edges corresponding to the litterals in $W_j$. Hence if we set $x_i$ to be 'True' precisely when $P$ avoids all vertices of the kind $u_{i,j}$ and false if it avoids all vertices of the kind $\bar{u}_{i,j}$ then we obtain a satisfying truth assignment for $\cal F$. Conversely, if $t$ is a satisfying truth assignment for $\cal F$ then we obtain the desired good completion $D'$ of $H'$ as follows:

\begin{itemize}
\item  For each $i\in [n]$: If $t(x_i)='False'$ then orient each edge $u_{i,j}v_{i,j}$, $j\in [p_i]$ from $u_{i,j}$ to $v_{i,j}$ and each edge $\bar{u}_{i,j}\bar{v}_{i,j}$, $j\in [q_i]$ from $\bar{v}_{i,j}$ to $\bar{u}_{i,j}$. If $t(x_i)='true'$, then orient oppositely to the above.
\item For each $j\in [m]$: If the litteral $\ell_{j,h}$ is true according to $t$, then orient the edge $a_{j,h}b_{j,h}$ from $a_{j,h}$ to $b_{j,h}$ and otherwise orient is from $b_{j,h}$ to $a_{j,h}$.
\end{itemize}

It is easy to check that is acyclic and contains an $(s,t)$-path, \qed

The mixed graph that we construct will never have a $(t,s)$-path as all arcs are directed into $t$. So this implies:

\begin{corollary}
It is NP-complete to decide whether a mixed graph  can be completed to be acyclic and have either an $(s,t)$-path or a $(t,s)$-path.
\end{corollary}

Even if the input acyclic mixed graph has both a (mixed) $(s,t)$-path and a (mixed) $(t,s)$-path the problem above  is still NP-complete. To see this take two copies $M_1,M_2$ of the mixed graph $H'$ from the proof of Theorem \ref{thm:mixedpath} and denote the copies of $s,t$ in these $s_1,t_1$, respectively $s_2,t_2$. Now add two new vertices $s^*,t^*$, and the arcs $s^*s_1,t_1t^*,t^*s_2,t_2s^*$. It is easy to check that the new mixed graph $H''$ is acyclic and has a completion to an acyclic digraph with an $(s^*,t^*)$-path or a $(t^*,s^*)$-path if and only if $H'$ has an acyclic orientation with an $(s,t)$-path.

%What happens if we require that the input satisfies that $\stackrel{\leftrightarrow}{P}$ has both an $(s,t)$-path and a $(t,s)$-path?

\section{Non-separating branchings in acyclic digraphs}\label{sec:nonsepBr}

A subdigraph  $D'=(V,A')$ of a digraph $D=(V,A)$ is {\bf non-separating} if the digraph $D^*=(V,A\setminus A')$ is connected.
It is well-known that one can decide in polynomial time whether an undirected graph $G$ has a pair of edge-disjoint spanning trees. This follows from matroid theory or the short direct proof in \cite{kaiserDM312} of Tutte's tree-packing theorem. By Edmonds's\ branching theorem \cite{edmonds1973} one can also decide in polynomial time whether a digraph contains two arc-disjoint out-branchings. In a personal communication Thomass\'e asked the following problem around 2008.
\begin{problem}
What is the complexity of deciding whether a digraph has a non-separating out-branching?
\end{problem}

This was answered by the authors of the present paper as follows.
\begin{theorem}\cite{bangTCS438}
  It is NP-complete to decide whether a digraph $D=(V,A)$ has a non-separating out-branching.
\end{theorem}

We now show that Thomass\'e's problem can be solved in polynomial time when $D$ is acyclic.

\begin{theorem}
  \label{polnonsepoutbranchingDAG}
  There exists a polynomial algorithm $\cal A$ which given an acyclic digraph $D=(V,A)$ either finds a non-separating out-branching or a certificate that $D$ has no such out-branching.
\end{theorem}

\pf
First observe that if $D$ has more than one vertex of in-degree zero then it has no out-branching and $\cal A$ can just out-put a certificate consisting of two vertices of in-degree zero. Hence we can assume that $D$ has a unique vertex $s$ of in-degree zero. Now let ${\cal O}=s=v_1,v_2,\ldots{},v_n$ be an arbitrary acyclic ordering of $D$ and note that every vertex $v_i$ with $i\geq 2$ has all its in-neighbours before it in $\cal O$. Furthermore, if we fix one arc entering each $v_i$, $i\geq 2$, then these arcs form an out-branching from $s$ in $D$.
This implies that $D$ has a non-separating out-branching if and only if $UG(D)$ has a spanning tree which uses at most $d^-(v_i)-1$ of the arcs entering $v_i$ for each $i\geq 2$. Let $M_1=(A,{\cal I}_1)$ be the matroid where a set \AY{$X$} of arcs is independent if and only if $X$ contains no more than $d^-(v_i)-1$ of the arcs entering $v_i$ for each $i\geq 2$ and let $M_2$ be the circuit matroid on $UG(D)$ (a set $Y$ of arcs is independent if and only if there is no cycle of $UG(D)$ all of whose arcs are contained in $Y$). Then it is easy to check that a set $Z$ of arcs is independent in both $M_1$ and $M_2$ if and only if $D-Z$ has an out-branching from $s$ and $Z$ induces a forrest in $UG(D)$. Hence $D$ has a non-separating out-branching if and only if the maximum size of a common independent set of $M_1$ and $M_2$ has size $n-1$. This can be checked by any polynomial algorithm for matroid intersection and in the case when the answer is no that algorithm can also deliver a certificate (based on Edmonds's matroid intersection theorem \cite{edmonds1970a}) for the non-existence of a non-separating out-branching.
\qed

By Theorem \ref{polnonsepoutbranchingDAG} and the easy fact that an undirected graph $G$  has an acyclic orientation with a non-separating out-branching if and only if $G$ has a pair of edge-disjoint spanning trees the following problem is natural to ask.

\begin{problem}
What is the complexity of deciding whether a mixed graph $P$ has an acyclic completion which contains a non-separating out-branching?
\end{problem}

We also know that deciding the existence of an out-branching $B^+_s$ which is arc-disjoint from some in-branching $B^-_t$ is polynomial for acyclic digraphs \cite{bangjGT42,bercziIPL23}. Hence it is tempting to ask about 
the complexity of deciding whether a mixed graph $M$ has an acyclic completion
which contains a pair $B^+_s,B^-_t$ of \AY{respectively} an out- and an in-branching.
However, already the case when $M$ has no arcs (is an undirected graph) seems very challenging and even the case when $M$ is the union of two edge-disjoint spanning trees is open and the problem is closely related to so-called rigidity matroids, see \cite{bangJGT96,bangJGT100}.\\

{\noindent}{\bf Acknowledgements}: J\o{}rgen Bang-Jensen thanks Jing Huang and Matthias Kriesell for stimulating discussions regarding the topics of Sections \ref{sec:outbr} and \ref{sec:stpath}.

%\bibliography{refs}

\end{document}